\documentstyle[a4,12pt]{article}
\begin{document}
\newtheorem{theorem}{Theorem}
\newtheorem{lemma}{Lemma}
\newtheorem{corollary}{Corollary}
\newcommand{\cl}{\mathop{\rm cl}\nolimits}
\hyphenation{semi-dihedral}
\hyphenation{gene-rated}
\hyphenation{gene-rators}
\hyphenation{pro-duct}
\hyphenation{pro-ducts}
\hyphenation{multi-pli-ca-ti-on}

\title{Wreath Products in the Unit Group of Modular Group
       Algebras of 2-groups of Maximal Class}
\author{Alexander B. Konovalov}
\date{}
\maketitle

\begin{abstract}
We study the unit group of the modular group algebra $KG$, where
$G$ is a \mbox{2-group} of maximal class.

We prove that the unit group of $KG$ possesses a section isomorphic
to the wreath product of a group of order two with the commutator
subgroup of the group $G$.

MSC2000: Primary 16S34, 20C05; Secondary 16U60

Keywords: Group algebras, unit groups, nilpotency class, wreath
products, 2-groups of maximal class

\end{abstract}

\section{Introduction}

Let $p$ be a prime number, $G$ be a finite $p$-group and $K$ be a
field of characteristic \nolinebreak $p$. Denote by
$\Delta=\Delta_{K}(G)$ the augmentation ideal of the modular group
algebra $KG$. The group of normalized units $U(G)=U(KG)$ consists of
all elements of the type $1+x$, where $x \in \Delta$. Our further
notation follows \cite{Sh1}.

Define Lie-powers $KG^{[n]}$ and $KG^{(n)}$ in $KG$: $KG^{[n]}$ is
two-sided ideal, generated by all (left-normed) Lie-products
$[x_{1},x_{2}, \cdots, x_{n}], x_{i} \in KG$, and $KG^{(n)}$ is
defined inductively: $KG^{(1)}=KG, KG^{(n+1)}$ is the associative
ideal generated by $[KG^{(n)},KG]$.  Clearly, for every
$n \quad KG^{(n)} \supseteq KG^{[n]}$, but equality need not hold.

For modular group algebras of finite $p$-groups $KG^{(|G'|+1)}=0$ \,
\cite{Sh3}. Then in our case finite lower and upper Lie nilpotency
indices are defined:
$$ t_{L}(G)=min\{n:KG^{[n]}=0\}, \qquad
   t^{L}(G)=min\{n:KG^{(n)}=0\}.$$

It is known that $ t_{L}(G) = t^{L}(G) $ for group algebras over the
field of characteristic zero \cite{PPS}, and for the case of
characteristic $p>3$ their coincidence was proved by A.~Bhandari and
I.~B.~S.~Passi \cite{Bhandari-Passi}.

Consider the following normal series in $U(G)$:
$$U(G)=1+\Delta \supseteq 1+\Delta(G\,') \supseteq
1+\Delta^{2}(G\,') \supseteq 1+\Delta^{t(G\,')}(G\,'),$$
where $t(G\,')$ is the nilpotency index of the augmentation
ideal of $KG\,'$.

An obvious question is whether does exist a refinement for this
normal series. There were two conjectures relevant to the question
above.

The first one, as it was stated in \cite{Sh1}, is attributed to
A.~A.~Bovdi and consists in the equality $\mbox{cl}\,U(G)=t(G\,')$,
i.e. this normal series doesn't have a refinement.
In particular, C.~Baginski \cite{Ba} proved that $\mbox{cl}\,U(G)=p$
if $|G\,'|=p$ (in case of cyclic commutator subgroup
$t(G\,')=|G\,'|$).  A.~Mann and A.~Shalev proved that $\cl U(G) \le
t(G\,')$ for groups of class two \cite{Mann-Sh}.

The second conjecture was suggested by S.~A.~Jennings \cite{Je3} in
a more general context, and in our case it means that
$\mbox{cl}\,U(G)=t_{L}(G)-1$. Here N.~Gupta and F.~Levin proved
inequality $\le$ in \cite{Gupta-Levin}.

The first conjecture was more attractive and challenging, since
methods for the systematic computation of the nilpotency index of the
augmentation ideal $t(G')$ were more known than such ones for the
calculation of the lower Lie nilpotency index (for key facts see, for
example, \cite{Je1}, \cite{Koshita}, \cite{Motose-Ninom}, \cite{Sh5},
\cite{Sh3}).

Moreover, A.~Shalev \cite{Sh1} proved that these two
conjectures are incompatible in general case, although
$t(G\,')=t_{L}(G)-1$ for some particular families of groups,
including 2-groups of maximal class. Later using computer Coleman
managed to find counterexample to Bovdi's conjecture (cf.
\cite{Sh4}), and the final effort in this direction was made bu X.~Du
\cite{Du} in his proof of Jennings conjecture.

Study of the structure of the unit group of group algebra and its
nilpotency class raised a number of questions of independent
interest, in particular, about involving of different types of wreath
products in the unit group (as a subgroup or as a section).

In \cite{Col-Pas} D.~Coleman and D.~Passman proved that for
non-abelian finite $p$-group $G$ a wreath product of two groups of
order $p$ is involved into $U(KG)$. Later this result was generalized
by A.~Bovdi in \cite{Bovdi3}. Among other related results it is worth
to mention \cite{Mann-Sh}, \cite{Mann}, \cite{Sh6}.

It is also an interesting question whether $U(KG)$ posesses a given wreath
product as a subgroup or only as a section, i.e. as a factor-group of a 
certain subgroup of $U(KG)$. Baginski in \cite{Ba} described all $p$-groups,
for which $U(KG)$ does not contain a subgroup isomorphic to the wreath 
product of two groups of order $p$ for the case of odd $p$, and the case
of $p=2$ was investigated in \cite{Bovdi-Dokuchaev}.

The question whether $U(G)$ possesses a section isomorphic to the
wreath product of a cyclic group $C_{p}$ of order $p$ and the
commutator subgroup of $G$ was stated by A.~Shalev in \cite{Sh1}.
Since the nilpotency class of the wreath product $C_{p} \wr H$ is
equal to $t(H)$ - the nilpotency index of the augmentation ideal of
$KH$ \cite{Buckley}, this question was very useful for the
investigation of the first conjecture. In \cite{Sh2} positive answer
was given by A.~Shalev for the case of odd $p$ and a cyclic
commutator subgroup of $G$.

The present paper is aimed to extend the last result on 2-groups of
maximal class, proving that if $G$ is such a group then the unit
group of $KG$ possesses a section isomorphic to the wreath product of
a group of order two with the commutator subgroup of the group $G$.
We prove the following main result.

\begin{theorem} \label{theorem1}
  Let $K$ be a field of characteristic two, $G$
  be a 2-group of maximal class. Then the wreath product
  $C_{2} \wr G\,'$ of a cyclic group of order two and
  the commutator subgroup of $G$ is involved in $U(KG)$.
\end{theorem}

\section{Preliminaries}

We consider 2-groups of maximal class, namely, the dihedral,
semidihedral and generalized quaternion groups, which we denote
by $D_{n},S_{n}$ and $Q_{n}$ respectively. They are given
by following representations \cite{Ba-mga-of-2groups}:
\begin{eqnarray}
D_{n} &=& \langle a,b \,|\, a^{2^{n-1}}=1, b^{2}=1, b^{-1}ab=a^{-1} \rangle ,
\nonumber \\
S_{n} &=& \langle a,b \,|\, a^{2^{n-1}}=1, b^{2}=1, b^{-1}ab=a^{-1+2^{n-2}} \rangle ,
\nonumber \\
Q_{n} &=& \langle a,b \,|\, a^{2^{n-1}}=1, b^{2}=a^{2^{n-2}}, b^{-1}ab=a^{-1} \rangle ,
\nonumber
\end{eqnarray}
where $n \ge 3$ (We shall consider $D_{3}$ and $S_{3}$ as
identical groups).

We may assume that $K$ is a field of two elements, since in the case
of an arbitrary field of characteristic two we may consider its simple
subfield and corresponding subalgebra in $KG$, where $G$ is one of the 
groups $D_{n},S_{n}$ or $Q_{n}$.

Denote for
$ x = \sum_{g \in G} \alpha_{g} \cdot g, \, \alpha_{g} \in K $
by  $\mbox{Supp}\,x$ the set $\{ g \in G \,|\, \alpha_{g} \ne 0 \}$.
Since $K$ is a f\/ield of two elements,
$x \in 1+\Delta \Leftrightarrow |\mbox{Supp}\,x|=2k+1, k \in {\rm N}$.

Next, for every element $g$ in $G$ there exists unique representation
in the form $ a^{i}b^{j} $, where $ 0 \le i < 2^{n-1}, 0 \le j < 2 $.
Then for every $ x \in KG $ there exists unique representation in the
form $ x=x_{1} +x_{2}b $, where
$ x_{i} = a^{n_{1}}+ \dots + a^{n_{k_{i}}} $.
We shall call $ x_{1}, x_{2} $ {\it components} of $x$. Clearly,
$ x_{1} +x_{2}b = y_{1} +y_{2}b \Leftrightarrow x_{i} = y_{i}, i=1,2 $.

The mapping $ x \mapsto \bar x = b^{-1}xb $, which we shall call  
{\it conjugation}, is an automorphism of order 2 of the group algebra $KG$.
An element $z$ such that $z=\bar z$ will be called {\it self-conjugated}.

Using this notions, it is easy to obtain the rule of multiplication
of elements from $KG$, which is formulated in the next lemma.

\begin{lemma} \label{multiply}
  Let $ f_{1} + f_{2}b, h_{1} + h_{2}b \in KG $. Then
  $$ (f_{1} + f_{2}b)(h_{1} + h_{2}b) =
     (f_{1}h_{1} + f_{2} \bar h_{2} \alpha) +
     (f_{2} \bar h_{1} + f_{1}h_{2})b ,$$
  where $ \alpha=1 $ for $D_{n}$ and $S_{n}$, $ \alpha=b^{2} $ for $Q_{n}$.
\end{lemma}

We proceed with a pair of technical results.

\begin{lemma} \label{lemma1}
  An element $ z \in KG $ commute with $ b \in G $ if and only if
  $z$ is self-conjugated.
\end{lemma}

\begin{lemma} \label{lemma2}
  If $x$ and $y$ are self-conjugated, then $xy=yx$.
\end{lemma}

In the next lemma we find the inverse element for an element from
$U(KG)$.

\begin{lemma} \label{lemma3}
  Let $ f=f_{1} +f_{2}b \in U(KG) $. Then
  $ f^{-1} = ( \bar f_{1} + f_{2}b) R^{-1}$, where
  $R = f_{1} \bar f_{1} + f_{2} \bar f_{2} \alpha$, and
  $\alpha=1$ for $D_{n}, S_{n}$, $\alpha=b^{2}$ for $Q_{n}$.
\end{lemma}

{\bf Proof}.
Clearly, $R$ is a self-conjugated element of $ K \langle a \rangle $
of augmentation 1, hence $R$ is a central unit in $KG$. Then the lemma
follows since
$$     (f_{1} + f_{2}b)(\bar f_{1} + f_{2}b)=
  (\bar f_{1} + f_{2}b)      (f_{1} + f_{2}b) = R. $$

Now we formulate another technical lemma, which is easy to prove by
straightforward calculations using previous lemma.

\begin{lemma} \label{lemma4}
  Let $ f,h \in U(G), f=f_{1} + f_{2}b, h=h_{1} + h_{2}b $,
  and $ h= \bar h $ is self-conjugated. Let
  $R$ and $\alpha$ be as in the lemma \ref{lemma3}.
  Then $ f^{-1}hf = t_{1} +t_{2}b $, where
  $ t_{1} = h_{1} + h_{2} (f_{1}f_{2} +
                          \bar f_{1} \bar f_{2}) \alpha R^{-1} ,
          t_{2} = h_{2} ( \bar f_{1}^{2} + f_{2}^{2} \alpha) R^{-1} $.
\end{lemma}

Let us consider the mapping
$ \varphi(x_{1} +x_{2}b) = x_{1} \bar x_{1} + x_{2} \bar x_{2} \alpha $,
where $ \alpha $ was defined in the Lemma \ref{lemma4}. It is easy to
verify that such mapping is homomorphism from $U(KG)$ to
$ U(K \!  \langle a \rangle ) $ and, clearly, for every $x$ its image
$ \varphi(x) $ is self-conjugated.  Such mapping
$ \varphi : U(KG) \rightarrow U(K \! \langle a \rangle ) $ we will
call {\em norm}. We will also say that the norm of an element $x$ is
equal to $ \varphi(x) $.

\section{Dihedral and Semidihedral Group}

Now let $G$ be the dihedral or semidihedral group.
Note that in Lemma \ref{lemma4} the first component of
$f^{-1}hf$ is always self-conjugated. In general, the second one need
not have the same property, but it is self-conjugated in the case
when $ f_{1} + f_{2} = 1 $, where $ f = f_{1} + f_{2}b $. It is easy
to check that the set
$$H(KG) = \{ h_{1} + h_{2}b \in U(KG) \,|\, h_{1} + h_{2} = 1 \}$$
is a subgroup of $U(KG)$.

Now we define the mapping
$ \psi : H(KG) \rightarrow U(K \! \langle a \rangle ) $
as a restriction of $ \varphi(x) $ on $H(KG)$. For convenience we
also call it {\em norm}.

\begin{lemma} \label{lemma5}
  $ \ker \psi = C_{H}(b) = \{ h \in H(KG) \,|\, hb=bh \} $.
\end{lemma}

The proof follows from the Lemma \ref{lemma1} and the equality
$ \psi(h) = 1 + h_{1} + \bar h_{1} $ for $ h \in H(KG) $.

\begin{lemma} \label{lemma6}
  $ C_{H}(b) $ is elementary abelian group.
\end{lemma}

{\bf Proof}.
$ C_{H}(b) $ is abelian by the lemma \ref{lemma2}. Now,
$ h^{2} = f_{1}^{2} + f_{2}^{2} + ( f_{1}f_{2} + f_{1}f_{2} )b =
f_{1}^{2} + f_{2}^{2} = 1 + f_{1}^{2} + f_{1}^{2} = 1 $,
and we are done.

Note that for elements $ f, h \in C_{H}(b) $ we may obtain more simple
rule of their multiplication:
$ (f_{1} + f_{2}b) (h_{1} + h_{2}b) =
  (1 + f_{1} + h_{1}) + (f_{1} + h_{1}) b $.

Now we consider a subgroup in $G$ generated by $ b \in G $ and
$A = a+(1+a)b$. Note that  $ A \in H(KG) $. First we calculate the
norm of $ A : \psi(A) = 1 + a + \bar a $. Since $ a^{2^{n-1}}=1 $, we
have $ \psi(A)^{2^{n-2}}=1 $. Then $ A^{2^{n-2}} $ commute with $b$,
and $ A^{2^{n-1}}=1 $ by the Lemma \ref{lemma6}. Smaller powers of
$A$ have non-trivial norm, so they do not commute with $b$. Clearly,
order of $A$ is equal to $ 2^{n-2} $ or $ 2^{n-1} $. In the following
lemma we will show that actually only the second case is possible.

\begin{lemma} \label{ord-a}
  Let $ A = a+(1+a)b $. Then the order of $A$ is equal to $2^{n-1} $.
\end{lemma}

{\bf Proof}.
We will show that the case of $2^{n-2}$ is impossible since
$ \mbox{Supp}\,A^{2^{n-2}} \ne 1 $. We will use formula (8) from
\cite{Bovdi-Lakatos}, which describes $2^{k}$-th powers of an element
$x \in U(G)$:
$$ x^{2^{k}} = x_{1}^{2^{k}} + ( x_{2} \bar x_{2} )^{2^{k-1}} b^{2^{k}} + \sum \limits_{i=1}^{k-1}
(x_{2} \bar x_{2})^{2^{i-1}} (x_{1} + \bar x_{1})^{2^{k}-2^{i}} b^{2^{i}} + x_{2} (x_{1} + \bar x_{1})^{2^{k}-1} b. $$

Let us show that the second component of $ A^{2^{n-2}} $ is
non-trivial. Let us denote it by $ t_{2}(A^{2^{n-2}}) = t_{2} $. By
the cited above formula, $ t_{2} = (1+a)(a + \bar a)^{2^{n-2}-1} $,
where $ \bar a = a^{-1} $ for the dihedral group, and
$ \bar a = a^{-1+2^{n-2}} $ for the semidihedral group.

Now we consider the case of the dihedral group. We have
$$ t_{2} = (1+a)(a+a^{-1})^{2^{n-2}-1} =
           (1+a)(a^{-1}(1+a^{2}))^{2^{n-2}-1} =$$
$$ (1+a) a^{2^{n-2}+1} (1+a^{2})^{2^{n-2}-1} =
(a^{2^{n-2}+1} + a^{2^{n-2}+2}) (1+a^{2})^{2^{n-2}-1} .$$
Note that if $ X = \langle x \rangle , x^{2^{m}}=1 $, then
$ (1+x)^{2^{m}-1} = \sum \limits_{y \in X} y = \overline{X} $, where
for a set $X$ we denote by $ \overline{X} $ the sum of all its elements
\cite{Ba-mga-of-2groups}. Thus, we have
$$ (1+a^{2})^{2^{n-2}-1} = 1 + a^{2} + a^{4} + \dots + a^{2^{n-1}-2} =
\overline{\langle a^{2} \rangle } = \overline{G\,'} .$$
Then
$$ t_{2} = (a^{2^{n-2}+1} + a^{2^{n-2}+2}) \overline{\langle a^{2} \rangle } =
            a \overline{\langle a^{2} \rangle } +
              \overline{\langle a^{2} \rangle } =
              \overline{\langle a \rangle }, $$
and for the case of the dihedral group the lemma is proved.

Now we will consider the semidihedral group. We have
$$ t_{2} = (1+a)(a+a^{2^{n-2}-1})^{2^{n-2}-1} =
           (1+a)(a(1+a^{2^{n-2}-2}))^{2^{n-2}-1} = $$
$$ (a^{2^{n-2}-1} + a^{2^{n-2}}) (1+a^{2^{n-2}-2})^{2^{n-2}-1} ,$$
and the rest part of the proof is similar. Note that from
$ t_{2}(A^{2^{n-2}}) = \overline{\langle a \rangle } $ we can immediately
conclude that its first component is $ 1+\overline{\langle a \rangle } $,
since $A \in H(G) $.

To construct a section isomorphic to the desired wreath product,
first we take elements
$ b, b^{A}, b^{A^{2}}, \cdots, b^{A^{2^{n-2}-1}} $.
For every $k$ we have $ (b^{A^{k}})^{2} = 1 $.  By the Lemma
\ref{lemma4} all elements  $ b^{A^{k}} $ are self-conjugated, since
$A \in H(KG)$, and they commute each with other by the Lemma
\ref{lemma2}. So, we get the next lemma.

\begin{lemma} \label{lemma7}
  $ \langle  b, b^{A}, b^{A^{2}}, \cdots, b^{A^{2^{n-2}-1}} \rangle $
  is elementary abelian subgroup.
\end{lemma}

Now we can obtain elements $ b^{A^{k}} $, using the Lemma \ref{lemma4}.

\begin{lemma} \label{lemma8}
  Let $ A = a+(1+a)b , \quad R = \psi(A) = 1 + a + \bar a $.
  Then
  $$ b^{A^{k}} = 1 + R^{k} + R^{k}b, \qquad k=1,2,...,2^{n-2} .$$
\end{lemma}

{\bf Proof}.
First we obtain $ b^{A} $ by Lemma \ref{lemma4} with
$ h_{1}=0, h_{2}=1, f_{1}=a, f_{2}=1+a $.
We get
$$ b^{A} = (a(1+a) + \bar a(1+ \bar a)) R^{-1} +
           (\bar a^{2} + (1+a)^{2}) R^{-1} b = $$
$$ = (a + \bar a + a^{2} + \bar a^{2}) R^{-1} +
     (1 + a^{2} + \bar a^{2}) R^{-1} b =
     (R + R^{2}) R^{-1} + R^{2} R^{-1} b  =  1 + R + Rb .$$
Now let $ b^{A^{k}} = 1 + R^{k} + R^{k}b $. Using the same
method for $ h_{1}=1 + R^{k}, h_{2}=R^{k} $, we get
$ b^{A^{k+1}} = 1 + R^{k+1} + R^{k+1}b $, as required.

\begin{lemma} \label{lemma9}
  There exists following direct decomposition:
  $$ \langle b, b^{A}, b^{A^{2}}, \cdots, b^{A^{2^{n-2}-1}} \rangle  =
  \langle b \rangle  \times \langle b^{A} \rangle  \times \langle b^{A^{2}} \rangle  \times \cdots
  \times \langle b^{A^{2^{n-2}-1}} \rangle $$
\end{lemma}

{\bf Proof}.
We need to verify that the product of the form
$ b^{i_{0}} (b^{A})^{i_{1}} \cdots (b^{A^{k}})^{i_{k}} $,
where $k=2^{n-2}-1$, $i_m \in \{0,1\}$ and not all $ i_{m} $ are equal to zero,
is not equal to $ 1 \in G $. Clearly, multiplication by $b$ only
permute components. So, we may consider only products
without $b$ and proof that they are not equal to 1 or $b$.

Note that $ b^{A^{k}} $ are self-conjugated and lies in $H(KG)$.
From this follows the rule of their multiplication:
$$ (1 + R^{k} + R^{k}b)(1 + R^{m} + R^{m}b) =
    1 + R^{k} + R^{m} + (R^{k} + R^{m})b .$$
The product of more than two elements is calculated by the same
way:
$$ (b^{A})^{i_{1}} (b^{A^{2}})^{i_{2}} \cdots (b^{A^{k}})^{i_{k}} =
    1 + i_1 R + i_2 R^2 + \cdots + i_k R^k +
      ( i_1 R + i_2 R^2 + \cdots + i_k R^k )b .$$

Put $\gamma = i_1 R + i_2 R^2 + \cdots + i_k R^k $ and $ R = 1 + r $, 
where $ r = a + \bar a $. Then $\gamma$ could be written in the form 
$ \gamma = \mu + r^{j_{1}} + \cdots + r^{j_{k}} $, where $\mu \in \{0,1\}$ 
and $ j_{1} < j_{2} < \cdots < j_{k} = i_{k} $. Since 
$ (a + \bar a)^{2^{n-2}} = 0 $, $r$ is nilpotent and its smaller powers are
linearly independent, so $ r^{j_{1}} + \cdots + r^{j_{k}} \ne 0 $. From the
other side, it is easy to see that the support of $ r^{j_s} $ does not 
contain 1, so $ r^{j_{1}} + \cdots + r^{j_{k}} \ne 1 $. Hence 
$ \gamma \not \in \{0,1\} $, and the support of the product
$ (b^{A})^{i_{1}} (b^{A^{2}})^{i_{2}} \cdots (b^{A^{k}})^{i_{k}} $
contains elements different from 1 and $b$, which proves the lemma.

Now we are ready to finish the proof of Theorem \ref{theorem1} for the
dihedral and semidihedral groups. 
It was shown that $U(KG)$ contains the semi-direct product $F$ of
$ \langle b         \rangle \times
  \langle b^{A}     \rangle \times
  \langle b^{A^{2}} \rangle \times
  \cdots                    \times
  \langle b^{A^{2^{n-2}-1}} \rangle  $
and $ \langle  A \rangle  $. As was proved above,
the order of $A$ is $ 2^{n-1} $ and its $2^{n-2}$-th power commutes
with $b$. From this follows that the factorgroup
$ F / \langle A^{2^{n-2}} \rangle $ is isomorphic to $ C_{2} \wr G\,' $,
as required.

\section{Generalized Quaternion Group}

Now let $G$ be the generalized quaternion group. First we need to
calculate $ \mbox{cl}\,U(G) $. In fact, we need to know only
$ t_{L}(G) $, since $ \mbox{cl}\,U(G) = t_{L}(G) - 1 $ \, \cite{Du}.
Note that Theorem \ref{theorem2} is already known (see Theorem 4.3 in 
\cite{Bovdi-Kurdics}), but we provide an independent proof for the 
generalized quaternion group.

\begin{theorem} \label{theorem2}
  Let $G$ be the generalized quaternion group. Then
  $ \mbox{cl}\,U(G) = |G\,'| $.
\end{theorem}

{\bf Proof}.
First, $ \mbox{cl}\,U(G) \le |G\,'| $ by \cite{Sh-Bi}. Now we
prove that \mbox{$ t_{L}(G) \ge |G\,'| + 1 $}. To do this, we will
construct non-trivial Lie-product of the length $ 2^{n-2} = |G\,'| $.

Consider Lie-product $ [b, \underbrace{a, \cdots, a}_{k}] $
which we denote by $[b, k \cdot A]$. Clearly, $[b,a]=(a+a^{-1})b$,
and $(a+a^{-1})$ is central in $KG$. It is easy to prove by induction
that $[b, k \cdot a]=(a+a^{-1})^k b$, therefore the commutator
$$[b, (2^{n-2}-1) \cdot a] = a^{2^{n-2}+1}(1+a^2)^{2^{n-2}-1} b =
a^{2^{n-2}+1} \overline {\langle a^2 \rangle} b = 
a \overline {\langle a^2 \rangle} b $$ does not vanish.

From the Theorem \ref{theorem2} it follows that
$ t_{L}(G) = t^{L}(G) $ since
$$ \mbox{cl}\,U(G) = t_{L}(G) - 1 \le t^{L}(G) - 1 \le |G\,'| ,$$
confirming conjecture about equality of the lower and upper
Lie nilpotency indices (cf. \cite{Bhandari-Passi}).
From this we conclude that $G$ and $U(KG)$ have the same exponent,
using the theorem from \cite{Sh3} about coincidence of their exponents
in the case when $ t^{L}(G) \le 1+(p-1)p^{e-1} $, where 
$ p^{e}=\mbox{exp}\,G $ and $p$ is the characteristic of the field $K$. 
Note that these two statements regarding Lie nilpotency indices and 
exponent are also true for all 2-groups of maximal class. Using the 
technique described here we also may show that modular group algebras of
2-groups of maximal class are Lie centrally metabelian.

For a unit $A$ of $KG$ we denote by $(b, k \cdot A)$ the commutator
$(b, \underbrace{A, \cdots, A}_{k})$.
Now we need a pair of technical lemmas.

\begin{lemma} \label{lemma10}
  Let $ A \in U(KG), A^{2^{n-1}}=1, b \in G,
  b^{A^{i}} b^{A^{j}} = b^{A^{j}} b^{A^{i}} $
  for every $ i,j $, where $ b^{A^{i}} = A^{-i} b A^{i} $.
  Then for every $ k \in {\bf N} \quad (b, k \cdot A)^{2} = 1 $.
\end{lemma}

{\bf Proof}.
We use induction by $k$. By straightforward calculation,
$ (b,A)^{2} = 1 $. Now, let  $(b, k \cdot A)=X, X^{2}=1 $. Then
$ (X,A)=XX^{A} $. Since elements $ b^{A^{i}}, i \in {\bf N} $
commute each with other, $X$ and $ X^{A} $ also commute, and
$ (XX^{A})^{2} = 1 $.

\begin{lemma} \label{lemma11}
  Let $ A \in U(KG), A^{2^{n-1}}=1, b \in G,
  b^{A^{i}} b^{A^{j}} = b^{A^{j}} b^{A^{i}} $
  for every $ i,j $, where $ b^{A^{i}} = A^{-i} b A^{i} $.
  Then for every $k,m \in {\bf N} $
$$ (b, \underbrace{A, \cdots, A}_{k}, A^{2^{m}}) =
   (b, \underbrace{A, \cdots, A}_{k+2^{m}} ).$$
\end{lemma}

{\bf Proof}.
We use induction by $m$. First,
$ (b, k \cdot A,A)=(b,(k+1) \cdot A) $.
Let the statement holds for some $m$. Consider the commutator
$$ (b, k \cdot A, A^{2^{m+1}}) =
   (b, k \cdot A, A^{2^{m}})^{2}
   (b, k \cdot A, A^{2^{m}}, A^{2^{m}}) ,$$
since $ (x,yz)=(x,y)(x,z)(x,y,z) $. By the Lemma \ref{lemma10} the square
of the first commutator is 1, while the second is equal to
$ (b, (k+2^{m+1}) \cdot A) $.

This gives possibility to proof the next property of $U(KG)$.

\begin{lemma} \label{theorem3}
  Let $ A \in U(KG), A^{2^{n-1}}=1, b \in G,
  b^{A^{i}} b^{A^{j}} = b^{A^{j}} b^{A^{i}} $
  for every $ i,j $, where $ b^{A^{i}} = A^{-i} b A^{i} $.
  Then $ A^{2^{n-2}} $ commute with $b$.
\end{lemma}

{\bf Proof}.
We will show using induction by $m$ that the group commutator
$(b, A^{2^{m}}) = (b, \underbrace{A, \cdots, A}_{2^{m}})$,
so \mbox {$ (b,A^{2^{n-2}})=(b, \underbrace{A, \cdots, A}_{2^{n-2}}) = 1 $},
since $ \mbox{cl}\,U(G) = 2^{n-2} $.

First, $ (b,A^{2}) = (b,A)^{2} (b,A,A) = (b,A,A) $ by the Lemma
\ref{lemma10}.  Let $ (b,A^{2^{m}}) = (b, 2^{m} \cdot A) $. Then
$(b,A^{2^{m+1}}) = (b,A^{2^{m}})^{2} (b,A^{2^{m}},A^{2^{m}}) =
(b,A^{2^{m}},A^{2^{m}}) $ by the Lemma \ref{lemma10}. Using induction
hypothesis and Lemma \ref{lemma11}, we get
$$ (b, A^{2^{m}}, A^{2^{m}}) = (b, 2^{m} \cdot A, A^{2^{m}}) =
   (b, 2^{m+1} \cdot A).$$

Let us take an element $ A=a^{2^{n-3}+1} + (1+a)b $, where
$ a^{2^{n-1}}=1 $. Calculating $ A^{-1}hA $ for self-conjugated
$h$ by the Lemma \ref{lemma4}, we get a self-conjugated element again.
The norm of $A$ is $ \varphi(A) = 1 + a^{2^{n-2}+1} + a^{2^{n-2}-1}$,
so order of $ \varphi(A) $ is $2^{n-2}$, and from this we conclude that
the order of $A$ is great or equal to $2^{n-2}$. From the other side, it is
not greater then $2^{n-1}$, since $G$ and $U(KG)$  have the same exponent.
Moreover, if $A^{2^{n-2}} \ne 1 $, then $A^{2^{n-2}}$ commute with $b$
by lemma \ref{theorem3}, and it is necessary to know whether its lower
powers commute with $b$. As in the previous section, in the following lemma
we will exactly calculate the order of $A$.

\begin{lemma} \label{ord-a-q}
  Let $ A=a^{2^{n-3}+1} + (1+a)b $.
  Then the order of $A$ is equal to $ 2^{n-1} $.
\end{lemma}

{\bf Proof}.
The proof is similar to the proof of the lemma \ref{ord-a}. We will show
that $ \mbox{Supp}\,A^{2^{n-2}} \ne 1 $, calculating the second component
$ t_{2}(A^{2^{n-2}}) = t_{2} $. Using the same formula from
\cite{Bovdi-Lakatos}, we have:
$$ t_{2} = (1+a)(a^{2^{n-3}+1} + a^{-2^{n-3}-1})^{2^{n-2}-1} =
           (1+a)(a^{-2^{n-3}-1} (1+a^{2^{n-2}+2}))^{2^{n-2}-1} = $$
$$         (1+a)(a^{-2^{n-3}-1})^{2^{n-2}-1} (1+a^{2^{n-2}+2})^{2^{n-2}-1} =
           (1+a) a^{-2^{n-3}+1} (1+a^{2^{n-2}+2})^{2^{n-2}-1} = $$
$$         (a^{-2^{n-3}+1} + a^{-2^{n-3}+2}) (1+a^{2^{n-2}+2})^{2^{n-2}-1} .$$

Then,
$ (1+a^{2^{n-2}+2})^{2^{n-2}-1} =
  \overline{\langle a^{2} \rangle } =
  \overline{G\,'} $.
From this
$$ t_{2} = (a^{-2^{n-3}+1}+a^{-2^{n-3}+2}) \overline{\langle a^{2} \rangle }=
            a \overline{\langle a^{2} \rangle } +
              \overline{\langle a^{2} \rangle } =
              \overline{\langle a \rangle },$$
and the lemma is proved.

Now we calculate elements $ b^{A^{k}}, k=1,2, \cdots , 2^{n-2} $,
using Lemma \ref{lemma4}.

\begin{lemma} \label{lemma12}
  Let $ A = a^{2^{n-3}+1} + (1+a)b, \quad
  R = \varphi(A) = 1 + a^{2^{n-2}+1} + a^{2^{n-2}-1}$.
  Then
  $$ b^{A^{k}} = \beta \sum \limits_{i=-1}^{k-1} (b^{2}R)^{i} +
  (b^{2}R)^{k}b, \quad k=1,2,...,2^{n-2} ,$$
  where $ \beta = a^{2^{n-3}+1} + a^{-2^{n-3}-1} +
                  a^{2^{n-3}+2} + a^{-2^{n-3}-2} $.
\end{lemma}

{\bf Proof}.
Remember that for $Q_{n}$ in Lemma \ref{lemma4} $ f^{-1}hf = t_{1} +t_{2}b $,
where
$$ t_{1} = h_{1} + h_{2}(f_{1}f_{2} + \bar f_{1} \bar f_{2})b^{2}R^{-1},
\qquad  t_{2} = h_{2} ( \bar f_{1}^{2} + f_{2}^{2} b^{2}) R^{-1}.$$
First we obtain the second component. For $ A = f_{1} + f_{2}b $ we
have $ \bar f_{1}^{2} + f_{2}^{2} b^{2} = (a^{-2^{n-3}-1})^{2} +
(1+a)^{2} a^{2^{n-2}} = b^{2} (1+a^{2}+a^{-2}) = b^{2}R^{2} $.
Then the second component of $ b^{A} $ is $ b^{2}R^{2}R^{-1} = b^{2}R $.
Now it is easy to prove by induction that the second component of
$ b^{A^{k}} $ is $ (b^{2})^{k} R^{k} $. From this immediately follows
that $ A^{k}, k < 2^{n-2} $, doesn't commute with $b$, since
$ \mbox{ord}\,R = 2^{n-2} $.

Now we will calculate the first component. First, for the element $A$
expression of the form $ f_{1}f_{2} + \bar f_{1} \bar f_{2} $ is equal to
$ a^{2^{n-3}+1} (1+a) + a^{-2^{n-3}-1} (1+a^{-1}) =
  a^{2^{n-3}+1} + a^{-2^{n-3}-1} + a^{2^{n-3}+2} + a^{-2^{n-3}-2} $,
which we will denote by $ \beta $. Using the formula at the beginning of
the proof for $ h_{1}=0, h_{2}=1 $ we conclude that the first component of
$b^{A}$ is equal to $ \beta b^{2} R^{-1} $. Now let the first component
of $b^{A^{k}}$, where $ k < 2^{n-2}-1 $, is equal to
$ \beta \sum \limits_{i=-1}^{k-1} (b^{2}R)^{i} $. Taking into consideration
its previously calculated second component, we obtain that the first
component of $b^{A^{k+1}}$ is equal to
$$ \beta \sum \limits_{i=-1}^{k-1} (b^{2}R)^{i} +
   \beta (b^{2}R)^{k} b^{2} R^{-1} =
   \beta \sum \limits_{i=-1}^{k} (b^{2}R)^{i} .$$

Now we are ready to construct the subgroup, whose factorgroup is
isomorphic to the desired wreath product. Let us consider the subgroup
$ F_{1} $ in $ U(G) $ :
$$ F_{1} = \langle b, b^{A}, b^{A^{2}}, \cdots, b^{A^{2^{n-2}-1}} \rangle
           \langle A \rangle , $$
where $ b \in G, A = a^{2^{n-3}+1} + (1+a)b $, $A^{2^{n-2}}$ is the minimal
power of $A$ which commutes with $b$. Further, the subgroup
$ \langle b, b^{A}, b^{A^{2}}, \cdots, b^{A^{2^{n-2}-1}} \rangle $ is
abelian, and the intersection of subgroups $ \langle b \rangle ,
\langle b^{A} \rangle , \cdots, \langle b^{A^{2^{n-2}-1}} \rangle $
is $ \langle b^{2} \rangle  $. Moreover, the order of $A$ is $2^{n-1}$.

Let us take $F_{2}$ as a factorgroup of the group $F_{1}$ as follows:
$$ F_{2} = F_{1} / \langle b^{2} \rangle \langle A^{2^{n-2}} \rangle  .$$
It is clear, that $ \cl F_{2} \le 2^{n-2} = \cl U(G) $. If we will show
that actually we have equality $\cl F_{2} = 2^{n-2} = \cl (C_{2} \wr G\,')$,
then from this it will follow that $ F_{2} \cong C_{2} \wr G\,' $.

Let $ M = C_{2} \wr G\,' $ and $ \cl F_{2} = 2^{n-2} = \cl M $.
Let us assume that $ F_{2} \not \cong M $. Then there exists such
normal subgroup $ N \triangleleft M $, that $ M/N \cong F_{2} $,
since there exists a homomorphism $M \rightarrow F_{2}$, which is
induced by mapping of generators of $M$ into $F_{2}$. Since
$ |Z(M)|=2$, $ N \triangleleft M $, then $ N \cap Z(M) \ne \emptyset $,
so $ Z(M) \subseteq N $. In this case the nilpotency class
$ \cl F_{2} $ should be less then $ \cl M $, and we will get a 
contradiction.

To obtain the lower bound for the nilpotency class $ \cl F_{2} $
we will show that the commutator $(b, \underbrace{A \dots A}_{2^{n-2}-1})$
in $U(G)$ does not belong to the subgroup
$ \langle b^{2} \rangle \langle A^{2^{n-2}} \rangle $, so its image
in $F_{2}$ is nontrivial. By the lemma \ref{lemma11}
$ (b, \underbrace{A \dots A}_{2^{n-2}-1}) =
  (b, \underbrace{A \dots A}_{2^{n-3}-1}, A^{2^{n-3}}) = $
$ (b, \underbrace{A \dots A}_{2^{n-4}-1}, A^{2^{n-4}}, A^{2^{n-3}}) =
    \cdots =
  (b, A, A^{2}, A^{4}, \cdots, A^{2^{n-4}}, A^{2^{n-3}})$,
and we obtain more simple commutator of the length $n-1$. Further,
$ (b,A) = b^{-1}b^{A} = bb^{A} $, then 
$ (b,A,A^{2}) = b b^{A} b^{A^{2}} b^{A^{3}} $, and, by induction,
$$(b, A, A^{2}, \cdots, A^{2^{n-3}}) =
   b b^{A} b^{A^{2}} \cdots b^{A^{2^{n-2}-1}} =
   (bA^{-1})^{2^{n-2}} A^{2^{n-2}}.$$

It remains to show that $ (bA^{-1})^{2^{n-2}} $ does not contained
in $ \langle b^{2} \rangle \langle A^{2^{n-2}} \rangle $. Note that
$(bA^{-1})^{-1} = Ab^{3}$,   $(Ab^{3})^{2^{n-2}} = (Ab)^{2^{n-2}}$.

By the Lemma \ref{ord-a-q} the second component of $ A^{2^{n-2}} $ is
equal to $ \overline{\langle a \rangle } $. Note that it is not changed
under multiplication of $ A^{2^{n-2}} $ by $b^{2}$. The same method
could be used for calculation of $ (Ab)^{2^{n-2}} $. We have
$$ Ab = (1+a)b^{2} + a^{2^{n-3}+1}b =
       ( a^{2^{n-2}} + a^{2^{n-2}+1} ) + a^{2^{n-3}+1}b .$$
Then by the formula from \cite{Bovdi-Lakatos} the second component of
$(Ab)^{2^{n-2}}$ is equal to
$$ a^{2^{n-3}+1} ( a^{2^{n-2}+1} + a^{2^{n-2}-1} )^{2^{n-2}-1} =
   a^{2^{n-3}+1} ( a^{2^{n-2}-1} (1 + a^{2}) )^{2^{n-2}-1} = $$
$$ a^{2^{n-3}+1} ( a^{2^{n-2}-1} )^{2^{n-2}-1} (1 + a^{2})^{2^{n-2}-1} =
   a^{2^{n-3}+2} \overline{\langle a^{2} \rangle } =
                 \overline{\langle a^{2} \rangle } .$$
Thus, support of the second component of $ (Ab)^{2^{n-2}} $ does not
coincide with the support of the second component of  $ (A)^{2^{n-2}} $
and does not changes under multiplication of $(Ab)^{2^{n-2}}$ by $b^{2}$.
From this we conclude that
$(Ab)^{2^{n-2}} \not\in \langle b^{2} \rangle \langle A^{2^{n-2}} \rangle $.
This proves that the commutator $ (b, \underbrace{A \dots A}_{2^{n-2}-1}) $
also does not lies there. That is why
$ \cl F_{2} = 2^{n-2} = \cl (C_{2} \wr G\,') $, and
$ F_{2} \cong C_{2} \wr G\,' $, so the theorem is proved.

$$ \mbox{\bf Acknowledgements} $$
The research was supported by the Hungarian National Foundation for
Scientific Research Grant No.~T~025029. The author is grateful to 
Prof.~Ya.~P.~Sysak for drawing attention to the problem and helpful
suggestions, and to referee for his useful comments.

\vspace{2 cm}
Department of Mathematics and Economy Cybernetics, \par
Zaporozhye State University, Zaporozhye, Ukraine \par
\vspace{5pt}
Mailing address: \par
P.O.Box 1317, Central Post Office, Zaporozhye, 69000, Ukraine \par
E-mail : konovalov@member.ams.org

\end{document}